\input amstex
\documentstyle{amsppt}
\magnification=\magstep1

\loadbold
\NoBlackBoxes

\pageheight{8.125truein}
\pagewidth{5.5truein}
\hcorrection{0.375truein}
\vcorrection{0.5truein}

\input xy
\xyoption{matrix}\xyoption{arrow}
\def\edge{\ar@{-}}

\long\def\ignore#1{#1}

\def\Pbar{\overline{P}}
\def\aut{\operatorname{aut}}
\def\spec{\operatorname{spec}}

\def\O{{\Cal O}}

\def\Oq{{\Cal O}_q}
\def\Oqm{{\Cal O}_q (k^m)}
\def\Oqn{{\Cal O}_q (k^n)}
\def\Olam{{\Cal O}_{\boldsymbol \lambda}}

\def\Oqmn{{\Cal O}_q(M_{m,n}(k))}
\def\Oqmt{{\Cal O}_q(M_{m,t}(k))}
\def\Oqtn{{\Cal O}_q(M_{t,n}(k))}
\def\Oqtt{{\Cal O}_q(M_t(k))}

\def\OqGLt{{\Cal O}_q ( GL_t(k))}

\def\Mmn{M_{m,n}(k)}
\def\Mmt{M_{m,t}(k)}
\def\Mtn{M_{t,n}(k)}

\def\ovt{\overline{\theta}}

\def\H{{\Cal H}}

\def\A{{\Cal A}}

\def\I{{\Cal I}}

\def\KMurcia{{\bf 1}}
\def\GLC{{\bf 2}}
\def\GL{{\bf 3}}
\def\GLet{{\bf 4}}
\def\McRob{{\bf 5}}
\def\PaWa{{\bf 6}}
\def\Rig{{\bf 7}}

\topmatter

\title Prime ideals in certain quantum determinantal rings\endtitle

\author K. R. Goodearl and T. H. Lenagan \endauthor

\address Department of Mathematics, University of California, Santa
Barbara, CA 93106, USA\endaddress
\email goodearl\@math.ucsb.edu\endemail

\address Department of Mathematics, J.C.M.B., Kings Buildings, Mayfield
Road, Edinburgh EH9 3JZ, Scotland\endaddress
\email tom\@maths.ed.ac.uk\endemail

\thanks This research was partially supported by National Science
Foundation research grant DMS-9622876 and NATO Collaborative
Research Grant 960250.\endthanks

\abstract The ideal $\I_1$ generated by the $2\times 2$ quantum minors 
 in the coordinate algebra of
quantum matrices, $\Oq(M_{m,n}(k))$, is  investigated. 
Analogues of the First and Second Fundamental Theorems of Invariant
Theory are proved.  In particular, it is shown that $\I_1$ is a 
completely prime ideal, that is, $\Oq(M_{m,n}(k))/\I_1$ is an integral
domain, and that  $\Oq(M_{m,n}(k))/\I_1$ is the ring of coinvariants of a
coaction of $k[x,x^{-1}]$ on $\Oq
(k^m) \otimes \Oq (k^n)$, a tensor product of two  quantum
affine spaces.  There is a natural torus action on $\Oq(M_{m,n}(k))/\I_1$
induced by an $(m+n)$-torus action on $\Oq(M_{m,n}(k))$.  We identify the
invariant prime ideals for this action and deduce consequences for the
prime spectrum of $\Oq(M_{m,n}(k))/\I_1$.
\endabstract

\title Prime ideals in certain quantum determinantal rings \endtitle

\endtopmatter

\document

\head Introduction\endhead

Let $k$ be a field and let $q \in k^\times$.
The {\it coordinate ring of quantum $m\times n$ matrices},
 $\A:= \Oqmn$, is a deformation of the
classical coordinate ring of $m\times n$ matrices, $\O (M_{m,n}(k))$.  As
such it is a $k$-al\-ge\-bra generated by $mn$ indeterminates $X_{ij}$, for
$1\leq i\leq m$ and $1\leq 
j\leq n$,
subject to the relations
$$\aligned X_{ij}X_{lj} &= qX_{lj}X_{ij}\\
X_{ij}X_{is} &= qX_{is}X_{ij}\\
X_{is}X_{lj} &= X_{lj}X_{is}\\
X_{ij}X_{ls} - X_{ls}X_{ij} &= (q - q^{-1})X_{is} X_{lj}\endaligned
\qquad\aligned &\text{when\ } i<l;\\
&\text{when\ } j<s;\\
&\text{when\ } i<l\ \text{and\ }j<s;\\
&\text{when\ } i<l\ \text{and\ }j<s.\endaligned$$
In some references (e.g., \cite{\PaWa, \S3.5}), $q$ is replaced by
$q^{-1}$. When $q=1$ we recover 
$\O (M_{m,n}(k))$, which is 
the commutative polynomial algebra $k[X_{ij}]$. 

When $m=n$, the algebra $\A $ possesses a special element, the {\it quantum
determinant}, $D_q$, defined by 
$$D_q :=  \sum _{ \sigma \in S_n}
(-q)^{l(\sigma )}X_{1,\sigma (1)}
X_{2,\sigma (2)}\cdots
X_{n,\sigma (n)},$$
where $l(\sigma)$ denotes the number of inversions in the permutation
$\sigma$. The quantum determinant $D_q$ is a central element of $\A$ (see,
for example, \cite{\PaWa, Theorem 4.6.1}), and the localization
$\A[D_q^{-1}]$ is the {\it coordinate ring of the quantum general linear
group}, denoted $\Oq(GL_n(k))$.

If $I \subseteq \{ 1,\dots ,m\}$ and $J\subseteq \{ 1,\dots ,n\}$ with
$|I|=|J| =t $, let $D(I,J)$ denote the $t\times t$ quantum minor obtained
as the quantum determinant of the subalgebra of $\A$ obtained
by deleting generators $X_{ij}$ from the rows outside $I$ and from the
columns outside
$J$.  We write
$\I _t$ for the ideal generated by the $(t+1)\times (t+1)$ quantum minors
of
$\A$.  In \cite{\GL} it is proved that  $\A/\I _t$ is an
integral domain, for each $1\leq t \leq  \min\{m,n\}$.
Independently, Rigal \cite{\Rig} has shown that $\A/\I _1$ is a
domain; he also shows that $\A/\I _1$ is a maximal order in its
division ring of fractions.

There is an action of the torus $\H := (k^\times )^m \times (k^\times
)^n$ by $k$-al\-ge\-bra automorphisms on $\A$ such that 
$$(\alpha_1 , \dots ,\alpha_m , \beta_1 , \dots , \beta_n
)\cdot X_{ij} := \alpha_i \beta_j X_{ij}$$
 for all $i,j$.
The ideals $\I _t$ are easily seen to be invariant under $\H$; so there
is an induced action of $\H$ on the factor algebras $\A/\I _t$.  
In this paper, we study the prime ideal structure in the algebra 
$\A/\I
_1$, paying particular attention to the $\H$-in\-var\-i\-ant prime ideals.

\head 1. Complete primeness of $\I_1$\endhead

We give a direct derivation of the fact that $\A/\I _1$ is a
domain.  Although this is already established in both \cite{\GL} and
\cite{\Rig}, the proof we give here is so  much simpler and more transparent
than either of the previous proofs that we think it will be useful to
have it in a published form.

The {\it coordinate ring of 
quantum affine $n$-space}, denoted $\Oq (k^n)$,
 is defined to be the $k$-al\-ge\-bra generated by elements
$y_1, \dots , y_n$ subject to the relations $ y_iy_j = qy_j y_i$ for
each $1\leq i<j\leq n$.  It is well known that $\Oq (k^n)$ is an
iterated Ore extension, and thus, in particular, $\Oq (k^n)$ is a
domain.  Our strategy is to produce a homomorphism  of $\A$ into $\Oq
(k^m) \otimes \Oq (k^n)$.  This latter algebra can also be presented as
an iterated Ore extension and thus is a domain.  We show that $\I_1 $ is
the kernel of this map and so $\A/\I _1$ is a domain.

\proclaim{1.1. Theorem} The algebra $\Oqmn/\I_1 $ is isomorphic to a
subalgebra of the tensor product 
 $\Oq (k^m) \otimes \Oq (k^n)$.  In particular,  $\I_1 $ is
a completely prime ideal of  $\Oqmn$.
\endproclaim

\demo{Proof}
Let $\Oqm = k[y_1, \dots , y_m]$ and $\Oqn = k[z_1, \dots , z_n]$ be
the coordinate rings of quantum affine $m$-space and $n$-space,
respectively.   We
define  an algebra homomorphism $\theta : \A \rightarrow \Oq
(k^m) \otimes \Oq (k^n)$ such that $\theta (X_{ij}) = y_i\otimes z_j$ for
all $i,j$.  In order  that this does extend to a well-defined algebra
homomorphism, we must check that the elements $y_i\otimes z_j$ satisfy at
least the relations defining
$\A$.  These are routine verifications; for example, if $i<l$ and 
$j<s$ then 
$$(y_i \otimes z_j)(y_l \otimes z_s) = y_i y_l \otimes z_jz_s =  y_i y_l
\otimes qz_sz_j = q(y_i \otimes z_s)(y_l \otimes z_j),$$
while
$$(y_l \otimes z_s)(y_i \otimes z_j) = y_l y_i \otimes z_sz_j =  q^{-1}y_i y_l
\otimes z_sz_j = q^{-1} (y_i \otimes z_s)(y_l \otimes z_j).$$
Thus, 
$$(y_i \otimes z_j)(y_l \otimes z_s) -(y_l \otimes z_s)(y_i \otimes z_j)
= (q -q^{-1}) (y_i \otimes z_s)(y_l \otimes z_j),$$
so that the fourth relation of the introduction holds. One can also obtain
$\theta$ as the composition of the comultiplication $\A\rightarrow
\A\otimes\A$ with the tensor product of the quotient maps from $\A$ to
$\A/\langle X_{ij}\mid i>1\rangle$ and $\A/\langle X_{ij}\mid
j>1\rangle$. We shall pursue the latter point of view in the next section.

Thus, there exists a unique $k$-al\-ge\-bra homomorphism 
$$\theta : \A
\rightarrow \Oq
(k^m) \otimes \Oq (k^n)$$ such that $\theta (X_{ij}) = y_i\otimes z_j$
for all $i,j$.   If $i<l$ and
$j<s$ then the above calculations also show that $\theta ( X_{ij} X_{ls}
-q X_{is}X_{lj}) = 0$; thus $\I_1 \subseteq \ker (\theta )$.

Now, $\Oq (k^m) \otimes \Oq (k^n)$ is a domain, since it can be viewed
as a (multiparameter) quantum affine $(m+n)$-space with respect to the
generators
$y_1
\otimes 1,
\dots , y_m \otimes 1, 1\otimes z_1 , \dots , 1\otimes z_n$.  Hence,
$\ker (\theta ) $ is  a completely prime ideal of $ \A$.  We show
that $\I_1 = \ker (\theta )$, so that $\I_1$ is completely prime.  It
remains to show that the induced map $\ovt :\A /\I_1 \rightarrow \Oq
(k^m) \otimes \Oq (k^n)$ is injective.  
Let  ${\Cal S}$ denote the set of monomials $X_{i_{1}j_{1}}X_{i_{2}j_{2}}
\dots X_{i_{l}j_{l}}$ in $\A$ such that $i_1 \geq i_2 \geq \dots \geq i_l$
and $j_1 \leq j_2 \leq \dots \leq j_l$.  (We allow the monomial to be
equal to $1$ when $l=0$.)  We claim that the set $\overline{{\Cal S}}$ of
images  forms a spanning set of $\A /I_1$.  

	It suffices to show that an arbitrary monomial $C$ in $\A$ is
congruent modulo $\I_1$ to  a linear combination of monomials from ${\Cal
S} $.  We proceed by induction on the index sets, where row index
sequences $( i_1, i_2 , \dots , i_l )$ are ordered lexicographically
with respect to $\geq$, column index sequences $( j_1, j_2, \dots ,
j_l) $ are ordered lexicographically with respect to $\leq$, and pairs
of sequences are ordered lexicographically.  

If the claim fails, then it fails for a monomial 
$C= X_{i_{1}j_{1}}X_{i_{2}j_{2}}
\dots X_{i_{l}j_{l}}$ whose index set is minimal with respect to the
ordering given in the previous paragraph. In particular, $C\notin\Cal
S$.  Let $r$ be the first subindex such that either
$i_r < i_{r+1}$ or $j_r > j_{r+1}$.  

	If $i_r < i_{r+1}$ and  $j_r \geq  j_{r+1}$ then $C=\lambda C'$,
where $\lambda$ is either $1$ or $q$ and $C'$ is obtained from $C$ by
switching $X_{{i_r}j_r}$ and $X_{i_{r+1}j_{r+1}}$.  However, 
$$( i_1, 
, \dots , i_{r-1}, i_{r+1}, i_{r}, i_{r+2},       
, \dots , i_{l}) < ( i_1, i_2 , \dots , i_l )$$
in our ordering, so
$C'$ is congruent modulo $\I_1$ to a linear combination of elements of
${\Cal S}$. Then $C$ is congruent to such a linear combination,
contradicting our assumptions.  A similar contradiction occurs if
$i_r
\leq  i_{r+1}$  and  $j_r >  j_{r+1}$: this time, the row indices might not
change, but 
$$(j_1,
\dots , j_{r-1}, j_{r+1}, j_r, j_{r+2} , \dots , j_l) < (j_1, \dots ,
j_l),$$ 
so again we have a contradiction. Therefore, we must either have  $i_r <
i_{r+1}$ and  $j_r <  j_{r+1}$ or $i_r > i_{r+1}$    
and  $j_r >  j_{r+1}$.

	Suppose that $i_r < i_{r+1}$    
and  $j_r <  j_{r+1}$.  In this case, we have
$$X_{i_{r}j_{r}}X_{i_{r+1}j_{r+1}} -q
X_{i_{r+1}j_{r}}X_{i_{r}j_{r+1}} \in\I_1,$$ 
so that $C -qC' \in \I_1 $, where 
$$C' = X_{i_{1}j_{1}}\cdots X_{i_{r-1}j_{r-1}}
X_{i_{r+1}j_{r}}X_{i_{r}j_{r+1}}X_{i_{r+2}j_{r+2}} \cdots
X_{i_{l}j_{l}}.$$
  We obtain a contradiction as above.

	The final case is $i_r > i_{r+1}$ and $j_r >  j_{r+1}$, where we
have  $$X_{i_{r}j_{r}}X_{i_{r+1}j_{r+1}} -q^{-1}
X_{i_{r}j_{r+1}}X_{i_{r+1}j_{r}} \in\I_1.$$  Thus, $C -q^{-1}C' \in \I_1 $,
where  
$$C' = X_{i_{1}j_{1}}\cdots X_{i_{r-1}j_{r-1}}
X_{i_{r}j_{r+1}}X_{i_{r+1}j_{r}}X_{i_{r+2}j_{r+2}} \cdots
X_{i_{l}j_{l}},$$
 and once again we reach a contradiction.
 This finishes the proof of the claim and
establishes that $\overline{{\Cal S}}$ spans $\A /\I_1$.

Now, observe that in $\Oq (k^m) \otimes \Oq (k^n)$ we have 
$$\theta
(X_{i_{1}j_{1}}X_{i_{2}j_{2}}
\dots X_{i_{l}j_{l}}) = y_{i_{1}}y_{i_{2}} \dots y_{i_{l}} \otimes 
z_{j_{1}}z_{j_{2}}
\dots z_{j_{l}}.$$  The monomials $y_{i_{1}}y_{i_{2}} \dots y_{i_{l}}$
 with $i_1
\geq i_2 \geq \dots \geq i_l$ are linearly independent over $k$, and,
likewise, the monomials $z_{j_{1}}z_{j_{2}}
\dots z_{j_{l}}$ with $j_1 \leq j_2 \leq \dots \leq j_l$ are linearly
independent over $k$.  Hence, $\theta$ maps ${\Cal S}$ bijectively 
to a linearly
independent set in $\Oq (k^m) \otimes \Oq (k^n)$, so that
$\overline{{\Cal  
S}}$ is a linearly independent set in $\A/\I_1$.  Therefore, the map
$\ovt : \A/\I_1 \rightarrow \Oq (k^m) \otimes \Oq (k^n)$ maps the
k-basis $\overline{{\Cal
S}}$ bijectively onto a linearly independent set, so that $\ovt$ is
injective.
 \qed\enddemo 

\head 2. Coinvariants\endhead

 Theorem 1.1 has an invariant theoretic interpretation, which we discuss
in this section.  First, we outline what happens in the classical ($q=1$)
case.

\definition{2.1}
Let $M_{u,v}(k)$ denote the algebraic variety of $u\times v $ matrices
over
$k$. Fix positive integers $m,n$ and $t\le \min\{m,n\}$. The general linear
group
$GL_t(k)$ acts on
$\Mmt
\times
\Mtn$ via
        $$g\cdot (A,B):= (Ag^{-1}, gB).$$
Matrix multiplication yields a map
        $$\mu : \Mmt \times \Mtn \rightarrow \Mmn,$$
the image of which is the variety of $m\times n$ matrices with rank at most
$t$, and there is an induced map
$$\mu _* : \O (\Mmn) \rightarrow \O \bigl( \Mmt \times \Mtn \bigr) = \O
(\Mmt)
\otimes \O (\Mtn).$$
The First Fundamental Theorem of invariant theory identifies the fixed
ring of the coordinate ring $\O \bigl( \Mmt \times \Mtn \bigr)$ under the
induced action
of
$GL_t(k)$ as precisely the image of $\mu_*$. The Second Fundamental
Theorem
states that the kernel of $\mu _*$ is
$\I_t$, the ideal generated by the
 $(t+1)\times (t+1)$
minors of $\O (\Mmn)$, so that the coordinate ring of the
variety of $m\times n$ matrices of rank at most $t$ is $\O
(\Mmn)/\I_t$.
 As a consequence, since this variety is irreducible, the ideal
$\I_t$ is a prime ideal of $\O (\Mmn)$.   
\enddefinition

\definition{2.2}
We now proceed to explain the connection between Theorem 1.1 and the above
invariant theoretic point of view.  

The analog of $\mu_*$ is the $k$-al\-ge\-bra homomorphism
$$\theta_t : \Oqmn \rightarrow \Oqmt\otimes\Oqtn$$
induced from the comultiplication on $\Oqmn$, that is, 
$$\theta_t(X_{ij})=
\sum_{l=1}^t X_{il}\otimes X_{lj}$$
for $1\le i\le m$ and $1\le j\le n$.
The comultiplications on $\Oqmt$ and $\Oqtn$ yield $k$-al\-ge\-bra
homomorphisms
$$\align \rho_t &: \Oqmt \rightarrow \Oqmt\otimes\Oqtt\\
 &\hphantom{: \Oqmt} \hskip1.0truein \rightarrow \Oqmt\otimes \OqGLt\\
\lambda_t &: \Oqtn\rightarrow \Oqtt\otimes \Oqtn\\
 &\hphantom{: \Oqtn} \hskip1.0truein \rightarrow \OqGLt \otimes\Oqtn
\endalign$$ 
which make $\Oqmt$ into a right
$\OqGLt$-comodule and $\Oqtn$ into a left
$\OqGLt$-comodule. Since $\OqGLt$ is a Hopf algebra, the right comodule
$\Oqmt$ becomes a
\underbar{left} $\OqGLt$-comodule on composing $\rho_t$ with $1\otimes S$
followed by the flip (where $S$ denotes the antipode). Finally, the tensor
product of the two left
$\OqGLt$-comodules $\Oqmt$ and $\Oqtn$ becomes a left $\OqGLt$-comodule
via the multiplication map on $\OqGLt$. This comodule structure map,
$$\gamma_t : \Oqmt\otimes\Oqtn \rightarrow \OqGLt\otimes
\Oqmt\otimes\Oqtn,$$
can be described (using the Sweedler summation notation) as follows:
$$\gamma_t(a\otimes b)= \sum_{(a)} \sum_{(b)} S(a_1)b_{-1} \otimes
a_0\otimes b_0$$
where $\rho_t(a)= \sum_{(a)} a_0\otimes a_1$ and $\lambda_t(b)= \sum_{(b)}
b_{-1}\otimes b_0$ for $a\in \Oqmt$ and $b\in \Oqtn$. Note that for $t>1$,
the map $\gamma_t$ is not an algebra homomorphism, since neither the
antipode nor the multiplication map on $\OqGLt$ is an algebra
homomorphism. On the other hand, $\gamma_1$ is a $k$-al\-ge\-bra
homomorphism.

Recall that the {\it coinvariants} of the coaction $\gamma_t$ are the
elements $x$ in the tensor product $\Oqmt\otimes\Oqtn$ such that
$\gamma_t(x)=1\otimes x$. Quantum analogs of the First and Second
Fundamental Theorems would be the following:

{\it Conjecture 1.} The set of coinvariants of $\gamma_t$ equals the image
of $\theta_t$.

{\it Conjecture 2.} The kernel of $\theta_t$ is the ideal $\I_t$.

\noindent We have proved Conjecture 2 in \cite{\GL, Proposition 2.4}
(essentially; the cited result covers the case $m=n$, and the general case
follows easily by the method of \cite{\GL, Corollary 2.6}). However,
Conjecture 1 is open at present. Here we shall establish it in the case
$t=1$.
\enddefinition

\definition{2.3} Note that $\Oq(M_{m,1}(k))$ and $\Oq(M_{1,n}(k))$ are
quantum affine spaces on generators $X_{11}, X_{21},\dots, X_{m1}$ and
$X_{11}, X_{12}, \dots, X_{1n}$, respectively. In studying the case $t=1$,
it is convenient to replace $\Oq(M_{m,1}(k))$ and $\Oq(M_{1,n}(k))$ by
$\Oqm = k[y_1, \dots , y_m]$ 
 and $\Oqn
= k[z_1, \dots , z_n]$, respectively. Then $\theta_1$ becomes the
$k$-al\-ge\-bra homomorphism
$$\theta : \Oqmn\rightarrow \Oqm\otimes\Oqn, \qquad X_{ij}\mapsto
y_i\otimes z_j$$
used in the proof of Theorem 1.1. Next, the (quantum) coordinate ring of
$1\times1$ matrices is just a polynomial ring $k[x]$, and the (quantum)
coordinate ring of the $1\times1$ general linear group is the localization 
$k[x,x^{-1}]$. Thus, in the present case the coaction $\gamma_1$ becomes
the $k$-al\-ge\-bra homomorphism
$$\gather \gamma : \Oqm\otimes\Oqn \rightarrow k[x^{\pm1}]\otimes
\Oqm\otimes\Oqn,\\
y_i\otimes1\mapsto x^{-1}\otimes y_i\otimes1, \qquad 1\otimes
z_j\mapsto x\otimes1\otimes z_j. \endgather$$
\enddefinition

\proclaim{2.4. Theorem} 
The set of coinvariants of $\gamma$ is exactly the image of the algebra
$\Oqmn $ in
$\Oqm \otimes \Oqn$ under $\theta$.
\endproclaim

\demo{Proof} Clearly $\gamma\theta(X_{ij})= 1\otimes y_i\otimes z_j=
1\otimes \theta(X_{ij})$ for all $i,j$. Since $\theta$ and $\gamma$ are
$k$-al\-ge\-bra homomorphisms, it follows that the image of $\theta$ is
contained in the coinvariants of $\gamma$.

The algebra $\Oqm\otimes\Oqn$ has a basis consisting of pure tensors
$Y\otimes Z$ where $Y$ is an ordered monomial in the $y_i$ and $Z$ is an
ordered monomial in the $z_j$. Note that $\gamma(Y\otimes Z)=
x^{s-r}\otimes Y\otimes Z$ where $r$ and $s$ are the total degrees of $Y$
and $Z$, respectively. Hence, the images $\gamma(Y\otimes Z)$ are
$k$-linearly independent, and a linear combination $\sum_{l=1}^d \alpha_l
Y_l\otimes Z_l$ of distinct monomial tensors is a coinvariant for $\gamma$
if and only if each
$Y_l\otimes Z_l$ is a coinvariant.

Thus, we need only consider a single term
$$Y\otimes Z= y_{i_{1} }y_{i_{2 } } \cdots y_{i_r}\otimes
z_{ j_{1}}z_{ j_{2}}\cdots z_{ j_s}.$$
If $Y\otimes Z$ is a coinvariant, then because $\gamma(Y\otimes Z)=
x^{s-r}\otimes Y\otimes Z$ we must have $r=s$. Therefore
$$Y\otimes Z= \theta  
(X_{i_{1}j_{1}}X_{i_{2}j_{2}}\cdots X_{i_{r}j_{r}}),$$
which shows that $Y\otimes Z$ is in the image of $\theta$, as desired.
\qed\enddemo

\head 3. $\H$-in\-var\-i\-ant prime ideals of $\Oqmn/\I_1$ \endhead 

Under the mild assumption that our ground field $k$ is infinite, we
identify the
$\H$-in\-var\-i\-ant prime ideals of the domain
$\A/\I_1=
\Oqmn/\I _1$. (Recall that $\H$ denotes the torus $(k^\times)^m \times
(k^\times)^n$, acting on $\A$ as described in the introduction.) This
identifies the minimal elements in a stratification of
$\spec
\A/I_1$, and yields a description of this prime
spectrum as a finite disjoint union of commutative schemes corresponding
to Laurent polynomial rings.

\definition{3.1} Let $H$ be a group acting as automorphisms on a ring
$A$. We refer the reader to \cite{\KMurcia} for the definition of the {\it
$H$-stratification of $\spec A$}, and here recall only that the
$H$-stratum of $\spec A$ corresponding to an $H$-prime ideal $J$ is the set
$$\spec_JA := \{P\in\spec A \mid (P:H)=J\}.$$

In the case of the algebra $\A/\I_1$, we shall (assuming $k$ infinite)
identify the
$\H$-prime ideals -- they turn out to be the same as the $\H$-in\-var\-i\-ant
primes -- and thus pin down the minimum elements of the $\H$-strata.
Further, we shall see that each $\H$-stratum of $\spec\A/\I_1$ is
homeomorphic to the spectrum of a Laurent polynomial ring over an
algebraic extension of
$k$. This pattern is also known to hold for $\spec\A$ itself (at least
when $q$ is not a root of unity), but there the
$\H$-prime ideals have not yet been completely identified.
\enddefinition

\definition{3.2}  It turns out that if a generator $X_{ij}$ lies in an
$\H$-prime ideal $P$ of $\A$ containing $\I_1$, then either all the
generators from the same row, or all the generators from the same column
must also lie in
$P$.  This leads us to make the following definition.  

For subsets $I \subseteq \{ 1,\dots , m\}$ and $J\subseteq \{ 1,\dots
, n\}$, set  $$P(I,J):= \I_1 + \langle X_{ij} \mid i\in I\rangle + 
\langle X_{ij} \mid j\in J\rangle . $$
Obviously, $P(I,J)$ is an $\H$-in\-var\-i\-ant ideal of $\A$. We shall show that
$P(I,J)$ is (completely) prime, and hence $\H$-prime.
\enddefinition

\proclaim{Lemma} The factor algebra
$\A/P(I,J)$ is isomorphic to $\Oq (M_{m',n'}(k))/\I_1'$, where $m'=m-|I|$
and
$n'=n-|J|$, and  $\I_1'$ is the ideal generated by the $2\times 2$
quantum minors of $\Oq (M_{m',n'}(k))$. 
Hence, $P(I,J)$ is a completely prime ideal of $\A$.
\endproclaim 
 
\demo{Proof}
The second statement follows immediately from the first statement and
Theorem 1.1.

Set $I' := \{1, \dots , m\}\backslash I$, and $J' := \{1, \dots ,
n\}\backslash J$, and let $\A'$ be the $k$-subalgebra of $\A$ generated by
the $X_{ij}$ for $i\in I'$ and $j\in J'$. Note that $\A'$ is
isomorphic to
$\Oq (M_{m',n'}(k))$.  Let $\I_1'$ be the ideal of $\A'$ generated by the
$2\times 2$ quantum minors of $\A'$; that is, those for which both row
indices are in $I'$ and both column indices are in $J'$.  Obviously, $\I_1'
\subseteq \A'
\cap \I_1$, so that the inclusion $\A' \hookrightarrow \A$ induces a
$k$-al\-ge\-bra homomorphism $f: \A' /\I_1' \rightarrow \A / P(I,J)$.  It
suffices to show that $f$ is an isomorphism.  

The factor $\A / P(I,J)$ is generated by the cosets of those $X_{ij}$
with $i\in I'$ and $j\in J'$, since $X_{ij} \in P(I,J)$ whenever $i\in I$
or $j\in J$.  These cosets are all in the image of $f$; so $f$ is
surjective.  

Observe that there exists a $k$-al\-ge\-bra homomorphism $g: \A \rightarrow
\A'$ such that $g(X_{ij}) = X_{ij}$ when $i\in I'$ and $j\in J'$, and
$g(X_{ij}) =0 $ otherwise.  To see this, note that the only problematic
relations are those of the form
$X_{ij}X_{ls}-X_{ls}X_{ij}=(q-q^{-1})X_{is}X_{lj} $ for $i<l$ and $j<s$.
However, if $i\not\in I' $ then $X_{ij}$ and $X_{is}$ both map to zero,
and the relation maps to $0=0$.  Likewise, this happens in all cases
except when $i, l \in I'$ and $j,s\in J'$:  in this case, the relation
above maps to a relation in $\A '$.

Consider  a $2\times 2$ quantum minor in $\A$ of the form
$D= X_{ij}X_{ls}-qX_{is}X_{lj}$ where $i<l$ and $j<s$.  If $i\not\in I'$
then both $X_{ij}$ and $X_{is}$ are in $\ker (g)$ , so that
$D \in \ker (g)$.  Likewise,
$g(D) =0$ when $l\not\in I'$, or $j\not\in J'$,
or $s\not\in J'$.  On the other hand, $g(D)=D$ when $i,l\in I'$ and
$j,s\in J'$. Further,
$g(X_{ij} ) =0$ when
$i\in I$ or
$ j\in J$.  Hence, $g(P(I,J)) \subseteq \I_1'$.

Therefore, $g$ induces a $k$-al\-ge\-bra homomorphism $\overline{g}: \A /
P(I,J) \rightarrow \A ' / \I_1'$.  Both of these algebras are generated
by the cosets corresponding to those $X_{ij}$ such that $i\in I'$ and
$j\in J'$.  It follows that both $f\overline{g}$ and $\overline{g}f$ are
identity maps, since both $f$ and $\overline{g}$ preserve these cosets. 
Hence, $f$ is an isomorphism.
 \qed\enddemo

Somewhat suprisingly, the $P(I,J)$ turn out to be the only
$\H$-prime ideals of $\A$ that contain $\I_1$.  The following
lemma will be helpful in establishing this fact.

\proclaim{3.3. Lemma} Let $i,s\in \{1,\dots,m\}$ and $j,t\in
\{1,\dots,n\}$. Then there exist scalars $\alpha\in \{1,q^{\pm1},
q^{\pm2}\}$ and $\beta\in \{1,q^{\pm1}\}$ such that $X_{ij}X_{st}- \alpha
X_{st}X_{ij}$ and $X_{ij}X_{st}-\beta X_{it}X_{sj}$ lie in
$\I_1$. In particular, the cosets $X_{ij}+\I_1$ are all normal
elements of $\A/\I_1$.\endproclaim

\demo{Proof} If $i=s$, then in view of the relations in $\A$ we can take
$\alpha=\beta$ to be $q$, 1, or $q^{-1}$ (depending on whether
$j<t$ or
$j=t$ or $j>t$). Similarly, if $j=t$, we can take $\alpha \in
\{1,q^{\pm1}\}$ and $\beta=1$.

If $i<s$ and $j>t$, or if $i>s$ and $j<t$, then $X_{ij}$ and $X_{st}$
commute, and we can take $\alpha=1$. On the other hand, one of
$X_{it}X_{sj}- q^{\pm1} X_{ij}X_{st}$ is a $2\times2$ quantum minor, and
so we can take $\beta$ to be $q$ or $q^{-1}$.

Now suppose that $i<s$ and $j<t$. Then $X_{ij}X_{st}- q
X_{it}X_{sj}$ is a quantum minor, and we can take $\beta=q$. But
$X_{st}X_{ij}- q^{-1}X_{it}X_{sj}$ is also a quantum minor, so we have
$X_{ij}X_{st} \equiv qX_{it}X_{sj}\equiv q^2X_{st}X_{ij}$ $\pmod{\I_1}$,
and hence we can take $\alpha=q^2$.

The remaining case follows from the previous one by exchanging $(i,j)$ and
$(s,t)$, and then the final statement of the lemma is clear. 
\qed\enddemo

\proclaim{3.4. Proposition}  Assume that $k$ is an infinite field.  Then
the $\H$-prime ideals of $\Oqmn$ that contain $\I_1$ are precisely the
ideals
$P(I,J)$.
\endproclaim  

\demo{Proof}
By Lemma 3.2, we know that the ideals $ P(I,J)$ are
$\H$-prime.  Consider an arbitrary $\H$-prime ideal $P$ of $\A$ that
contains
$\I_1$. If all of the $X_{ij}$ are in $P$ then $P$ must be the maximal
ideal generated by the $X_{ij}$.  In that case, $P = P(I,J)$, where $I=\{1,
\dots , m\}$ and $J=\{ 1, \dots ,n\}$.   Hence, we may assume that not
all $X_{ij}$ are in $P$.
Set 
$$\align I &= \{ i\in \{1,\dots,m\} \mid X_{ij} \in P \;{\text{for all}}\;
j\}\\
J &= \{ j\in \{1,\dots,n\} \mid X_{ij} \in P \;{\text{ for all}}\; i\}.
\endalign$$
  
We first show that $X_{ij} \in P$ if and only if $i\in I $ or
$j\in J$.  Certainly, if $i\in I $ or $j\in J$ then $X_{ij} \in P$, by
the definition of $I$ and $J$.  Suppose that there exists an $X_{ij} \in
P$ such that $i\not\in I$ and $j\not\in J$.  Then there exists an
index $s\neq i$ such that $X_{sj}\not\in P$ and also there exists an
index $t\neq j$ such that $X_{it} \not\in P$.  By Lemma 3.3, there is a
nonzero scalar $\beta\in k$ such that
$X_{ij}X_{st} - \beta X_{it} X_{sj}\in P$.  Thus, $X_{ij}
\in P$ would imply that $X_{it}X_{sj} \in P$.  However, $X_{it}$ and
$X_{sj}$ are $\H$-eigenvectors which, by Lemma 3.3, are normal modulo $P$.
Hence, because $P$ is $\H$-prime, $X_{it}X_{sj} \in P$ would imply
$X_{it}\in P$ or
$X_{sj}\in P$, contradicting the
choices of $s$ and $t$. Thus, we have established that $X_{ij} \in P$ if
and only if
$i\in I
$ or
$j\in J$.  Now $P(I,J) \subseteq P$, and we need to establish equality. 

Set $B:=
\A /P(I,J)$ and $\Pbar= P/P(I,J)$, and note that $B$ is a domain by Lemma
3.2.  
Write $Y_{ij}$ for the image of $X_{ij}$ in $B$.  The claim just
established implies that $Y_{ij} \notin \Pbar$ if $i\notin I$ or
$j\notin J$.
Recall from Lemma 3.3 that the $Y_{ij}$ 
scalar-commute among themselves.

        Now, $I\neq \{ 1, \dots , m\}$ and $J\neq \{ 1, \dots , n\}$,
since not all of the $X_{ij} $ are in $P$.  Let
$s\in \{1, \dots , m\}\setminus I$ and $t\in \{1, \dots ,
n\}\setminus J$ be minimal, and consider the localization $C:=
B[Y_{st}^{-1}]$.  Since $Y_{st} \not\in \Pbar$ there is an embedding of
$B$ into $C$, and $\Pbar C$ is an $\H$-prime ideal of $C$ such that $\Pbar
C\cap B= \Pbar$.

        Note that $Y_{ij} =0$ if $i<s$ or $j<t$.  If $i>s$ and $j>t$,
then we have $Y_{st}Y_{ij} - q Y_{sj}Y_{it} =0$, so that $Y_{ij} = q
Y_{st}^{-1} Y_{sj} Y_{it} $ in $C$.  Hence, $C$ is generated as an algebra
by
$Y_{st}^{\pm 1}$ together with $Y_{sj}$ for $j>t$ and $Y_{it}$ for $i>s$.
 Thus, $C$ is a homomorphic image of a localized multiparameter quantum
affine space
$\Olam(k^r)[z_1^{-1}]$, for $r=m-s+n-t +1$ and for a
suitable parameter matrix $\boldsymbol\lambda$.

        The standard action of the torus $\H_r:= (k^\times )^r$ on
$\Olam(k^r)$ has 1-dimensional eigenspaces generated by individual
monomials (here, we use the fact that $k$ is infinite).  Therefore, the
same holds for $C$. Hence, any nonzero $\H_r$-in\-var\-i\-ant ideal of $C$
contains a monomial, and so any nonzero $\H_r$-prime ideal of $C$ must
contain one of $Y_{s+1,t}, \dots, Y_{mt}$, $Y_{s,t+1}, \dots, Y_{sn}$.
Since $\Pbar C$ contains none of these elements, to show that $\Pbar C=0$
it suffices to establish that $\Pbar C$ is $\H_r$-prime. But $\Pbar C$ is
already $\H$-prime, so it is enough to see that the $\H_r$-in\-var\-i\-ant
ideals of $C$ are the same as the $\H$-in\-var\-i\-ant ideals. This will follow
from showing that the images of $\H$ and $\H_r$ in $\aut C$ coincide.

Since the $Y_{ij}$ are $\H$-eigenvectors, it is clear that the image of
$\H$ is contained in the image of $\H_r$. The reverse inclusion amounts to
the following statement:  

(*) Given any $
\alpha_{s}, \dots , \alpha_m , \beta_{t+1}, \dots , \beta_n 
\in 
k^\times$, there exists $h\in \H$ such that $h(Y_{it}) = \alpha_i Y_{it} $
for $i= s,
\dots , m$ and $h(Y_{sj}) = \beta_jY_{sj}$
for $j= t+1, \dots , n$.

Now, there exists $h_1 \in \H$ such that $h_1(X_{ij}) =
X_{ij}$ for all  $i,j$ with  $i<s$, and $h_1 (X_{ij}) = \alpha_i X_{ij}$ 
for all $i,j $ with $i\geq s$.  Also, there exists $h_2 \in \H $ such
that $h_2 (X_{ij}) = X_{ij} $ for all $i,j$ with $j\leq t$ and $h_2
(X_{ij}) = \alpha_s^{-1} \beta_jX_{ij} $ for all $i,j$ with $j>t$.
Setting $h=h_1 h_2$ gives the desired element of $\H$, establishing (*).

	Therefore, $\Pbar C=0$, and so $\Pbar =0$.  This means that
$P=P(I,J)$. 
 \qed\enddemo

\proclaim{3.5. Corollary}  If the field $k$ is infinite, then $\Oqmn/\I_1$
has precisely
$(2^m -1)(2^n -1) +1$ distinct $\H$-prime ideals, all of which are
completely prime. Further, each $\H$-stratum of $\spec\Oqmn/\I_1$ is
homeomorphic to the prime spectrum of a Laurent polynomial ring over an
algebraic field extension of $k$.
\endproclaim

\demo{Proof} The first statement is clear from Proposition 3.4. The
second statement is not actually a corollary of the Proposition, but is
included to fill in the picture. It may be obtained from
\cite{\KMurcia, Theorems 5.3, 5.5} (all but the algebraicity of the
coefficient fields also follows from \cite{\GLet, Theorem 6.6}).
\qed\enddemo

\definition{3.6}
In particular, the corollary above
explains why in the algebra $\Oq (M_2(k))$ there
are precisely $10 = (2^2 -1)^2+1$ distinct $\H$-primes which contain
the quantum determinant. This fact was known previously by direct
enumeration of these primes. The remaining $\H$-primes correspond to
$\H$-primes of $\Oq(GL_2(k))$; there are 4 of these, as has long been
known. We can display the lattice of $\H$-prime ideals of $\Oq (M_2(k))$
as in the diagram below, where the symbols $\bullet$ and $\circ$ stand for
generators $X_{ij}$ which are or are not included in a given prime, and
$\square$ stands for the $2\times2$ quantum determinant. For example,
$\tbinom{\circ\bullet}{\bullet\circ}$ stands for the ideal $\langle X_{12},
X_{21}\rangle$, and
$(\square)$ stands for the ideal $\langle X_{11}X_{22} -qX_{12}X_{21}
\rangle$.

\goodbreak\midinsert
\ignore{
$$\xymatrixrowsep{2.4pc}\xymatrixcolsep{3.2pc}
\xymatrix{
 &&\tbinom{\bullet\bullet}{\bullet\bullet}\\
\tbinom{\bullet\bullet}{\circ\bullet} \edge[urr]
&\tbinom{\bullet\bullet}{\bullet\circ} \edge[ur]
&&\tbinom{\circ\bullet}{\bullet\bullet} \edge[ul]
&\tbinom{\bullet\circ}{\bullet\bullet} \edge[ull]\\
\tbinom{\bullet\bullet}{\circ\circ} \edge[u] \edge[ur] 
 &\tbinom{\circ\bullet}{\circ\bullet} \edge[ul] \edge[urr]
 &\tbinom{\circ\bullet}{\bullet\circ} \edge[ul] \edge[ur]
 &\tbinom{\bullet\circ}{\bullet\circ} \edge[ull] \edge[ur]
 &\tbinom{\circ\circ}{\bullet\bullet} \edge[ul] \edge[u]\\
 &\tbinom{\circ\bullet}{\circ\circ} \edge[ul] \edge[u] \edge[ur]
 &\left(\square\right) \edge[ull] \edge[ul] \edge[ur] \edge[urr]
 &\tbinom{\circ\circ}{\bullet\circ} \edge[ul] \edge[u] \edge[ur]\\
 &&\tbinom{\circ\circ}{\circ\circ} \save+<0pc,-2.5pc>
\drop{{\Cal H}{\operatorname{-spec}}\, {\Cal O}_q(M_2(k))} \restore 
\edge[ul] \edge[u] \edge[ur]
}$$
}
\endinsert

The corresponding $\H$-strata in $\spec \Oq(M_2(k))$ can be easily
calculated. For instance, if $q$ is not a root of unity, the strata
corresponding to $\tbinom{\circ\circ}{\circ\circ}$ and
$\tbinom{\circ\bullet}{\bullet\circ}$ are 2-dimensional, the strata
corresponding to $\tbinom{\bullet\bullet}{\circ\bullet}$,
$\tbinom{\bullet\bullet}{\bullet\circ}$,
$\tbinom{\circ\bullet}{\bullet\bullet}$, and
$\tbinom{\bullet\circ}{\bullet\bullet}$ are all 1-dimensional, and the
remaining 8 strata are singletons.
\enddefinition

\definition{3.7} We close with some remarks concerning catenarity.
(Recall that the prime spectrum of a ring $A$ is {\it catenary} provided
that for any comparable primes $P\subset Q$ in $\spec A$, all saturated
chains of primes from $P$ to $Q$ have the same length.) It is conjectured
that $\spec \Oqmn$ is catenary. In \cite{\GLC, Theorem 1.6}, we showed
that catenarity holds for any affine, noetherian, Auslander-Gorenstein,
Cohen-Macaulay algebra $A$ with finite Gelfand-Kirillov dimension,
provided $\spec A$ has normal separation. All hypotheses but  the last
 are known to hold for the algebra $\A= \Oqmn$. We can, at
least, say that the portion of $\spec\A$ above $\I_1$ -- that is, $\spec
\A/\I_1$ -- is catenary: In view of Lemma 3.3, $\A/\I_1$ is a homomorphic
image of a multiparameter quantum affine space $\Olam(k^{n^2})$, and $\spec
\Olam(k^{n^2})$ is catenary by \cite{\GLC, Theorem 2.6}.
\enddefinition

\enddocument